\documentclass[10pt]{article}
\usepackage{latexsym}
\usepackage{amsbsy}
\usepackage{amsfonts,mathrsfs,amssymb}
\usepackage{amsmath,amsthm}
\usepackage[latin1]{inputenc}
\usepackage{graphics,graphicx}
\usepackage{a4wide}

 \newtheorem{thm}{Theorem}[section]
 \newtheorem{defin}[thm]{Definition}
 \newtheorem{lem}[thm]{Lemma}
 \newtheorem{prop}[thm]{Proposition}
 \newtheorem{cor}[thm]{Corollary}
 \newtheorem{rem}[thm]{Remark}
 \newtheorem{ex}[thm]{Example}

 \newcommand{\bthm}{\begin{thm}}
 \newcommand{\ethm}{\end{thm}}
 \newcommand{\bd}{\begin{defin}}
 \newcommand{\ed}{\end{defin}}
 \newcommand{\blem}{\begin{lem}}
 \newcommand{\elem}{\end{lem}}
 \newcommand{\bcor}{\begin{cor}}
 \newcommand{\ecor}{\end{cor}}
 \newcommand{\bprop}{\begin{prop}}
 \newcommand{\eprop}{\end{prop}}
 \newcommand{\brem}{\begin{rem} \rm}
 \newcommand{\erem}{\end{rem}}
 \newcommand{\bex}{\begin{ex} \rm}
 \newcommand{\eex}{\end{ex}}

 \newcommand{\pr}{\noindent{\bf Proof. }}
 \newcommand{\ep}{\nolinebreak{\hspace*{\fill}$\Box$ \vspace*{0.25cm}}}

 \newcommand{\beq}{\begin{equation}}
 \newcommand{\eeq}{\end{equation} }

 \newcommand{\bea}{\begin{eqnarray}}
 \newcommand{\eea}{\end{eqnarray}}

 \newcommand{\beas}{\begin{eqnarray*}}
 \newcommand{\eeas}{\end{eqnarray*}}

 \newcommand{\beqs}{\begin{equation*}}
 \newcommand{\eeqs}{\end{equation*}}

 \newcommand{\bi}{\begin{itemize}}
 \newcommand{\ei}{\end{itemize}}

 \newcommand{\ben}{\begin{enumerate}}
 \newcommand{\een}{\end{enumerate}}

 \newcommand{\ba}{\begin{array}}
 \newcommand{\ea}{\end{array}}

 \newcommand{\R}{\mathbb R}

 \newcommand{\cD}{\ensuremath{{\cal D}}}
 
 \newcommand{\cF}{\ensuremath{{\cal F}}}

 \newcommand{\cL}{\ensuremath{{\cal L}}}

 \newcommand{\cS}{\ensuremath{{\cal S}}}

 \newcommand{\rep}{\ensuremath{\mbox{Re\,}}}
 \newcommand{\imp}{\ensuremath{\mbox{Im\,}}}

\begin{document}

 \title{
 Waves in fractional Zener type viscoelastic media
 }

 \author{
 Sanja Konjik
 \footnote{Faculty of Agriculture, Department of Agricultural Engineering, University of Novi Sad, Trg Dositeja Obradovi\' ca 8, 21000 Novi Sad, Serbia.
 Electronic mail: sanja\_konjik@uns.ac.rs}\\
 Ljubica Oparnica
 \footnote{Institute of Mathematics, Serbian Academy of Science and Art, Kneza Mihaila 35, 11000 Belgrade, Serbia.
 Electronic mail: ljubicans@sbb.co.rs}\\
 Du\v san Zorica
 \footnote{Faculty of Civil Engineering, University of Novi Sad, Kozara\v{c}ka 2a, 24000 Subotica, Serbia.
 Electronic mail: zorica@gf.uns.ac.rs}\\
 }

 \date{}
 \maketitle

 \begin{abstract}
 Classical wave equation is generalized for the case of
 viscoelastic materials obeying fractional Zener model instead of
 Hooke's law. Cauchy problem for such an equation is studied:
 existence and uniqueness of the fundamental solution is proven and
 solution is calculated.

 \vskip5pt

 \noindent {\bf Mathematics Subject Classification (2000):}
 Primary: 26A33; Secondary: 46F12, 74D05, 74J05

 \vskip5pt

 \noindent {\bf Keywords:}
 fractional derivatives, Laplace and Fourier
 transforms, fundamental solution, fractional Zener model
 \end{abstract}

 \section{Introduction} \label{sec:intro}

 The need for a use of derivatives of arbitrary real (or even complex) order,
 also called fractional derivatives, is permanently underlined by various
 applied problems. Although, the idea about differentiating with respect to
 noninteger order originates from the work of Leibnitz and Euler at the turn
 of the 18th century, the last 30 years have brought a real expansion of
 fractional calculus. Certainly, it is due to great number of possibilities
 for its application in many diverse fields such as mechanics,
 viscoelasticity, automatic control, signal processing, stochastic and
 finance, biomedical engineering, etc.
 For a detailed exposition of the theory of fractional calculus we refer to
 \cite{Podlubny} or \cite{SamkoKM}.

 In this paper we study existence and uniqueness of solutions for the system
 \begin{eqnarray}
 \frac{\partial}{\partial x}\sigma (x,t) &=&
 \frac{\partial^{2}}{\partial t^{2}}u(x,t),  \nonumber \\
 \sigma (x,t) + \tau {}_{0}D_{t}^{\alpha}\sigma (x,t)
 &=& \varepsilon (x,t) +{}_{0}D_{t}^{\alpha}\varepsilon (x,t),  \label{eq:system} \\
 \varepsilon (x,t) &=& \frac{\partial}{\partial x}u(x,t), \nonumber
 \end{eqnarray}
 where $x\in \mathbb{R}$, $t>0$, $\sigma$, $u$ and $\varepsilon$ are stress,
 displacement and
 strain, respectively, considered as functions of $x$ and $t$, $0<\tau<1$ is
 a constant and ${}_{0}D_{t}^{\alpha}$, $0\leq \alpha <1$, denotes the left
 Riemann-Liouville operator of fractional differentiation of order $\alpha$.
 This system describes waves in viscoelastic media.

 System (\ref{eq:system}) is subjected to initial conditions
 \begin{equation} \label{eq:ic}
 \begin{array}{c}
 \displaystyle u(x,0) = u_{0}(x), \qquad \frac{\partial}{\partial t}
 u(x,0) = v_{0}(x), \\
 \displaystyle \sigma (x,0) =0, \qquad \varepsilon(x,0) =0,
 \end{array}
 \end{equation}
 as well as boundary conditions
 \begin{equation} \label{eq:bc}
 \lim_{x\to\pm \infty}u(x,t)=0, \qquad \lim_{x\to \pm \infty} \sigma(x,t)=0.
 \end{equation}
 As we shall show, system (\ref{eq:system}) can be reduced to the equation
 \begin{equation} \label{eq:fzwe}
 \frac{\partial^{2}}{\partial t^{2}}u(x,t) =L(t)
 \frac{\partial^{2}}{\partial x^{2}}u(x,t) ,\qquad x\in \mathbb{R}, t>0,
 \end{equation}
 which will be called wave equation for fractional Zener type viscoelastic
 media or shortly FZWE (Fractional Zener Wave Equation). Here $L$
 denotes a
 linear operator (of convolution type) acting on $\cS'(\mathbb{R})$, whose
 explicit form is given by
 $$
 L(t) =\cL^{-1} \left(\frac{1+s^{\alpha}}{1+\tau s^{\alpha}}\right) \ast_{t}.
 $$
 In fact, equation (\ref{eq:fzwe}) will be the subject of our consideration.
 System (\ref{eq:system}), and in particular (\ref{eq:fzwe}), generalize
 classical wave equation.

 Equation that describes waves occurring in the elastic media are obtained by
 using basic equations of elasticity (see \cite{AtanackovicGuran}). Our interest is
 restricted to an infinite elastic rod (one dimensional body), positioned
 along $x$ axis, that is not under influence of body forces. Equations of
 elasticity then read%
 \begin{eqnarray}
 \frac{\partial}{\partial x}\sigma(x,t) &=&
 \rho \frac{\partial^{2}}{\partial t^{2}}u(x,t) \label{eq:4} \\
 \sigma (x,t) &=& E\varepsilon (x,t) \label{eq:5} \\
 \varepsilon (x,t) &=& \frac{\partial}{\partial x}u(x,t) \label{eq:6}
 \end{eqnarray}
 where $x\in \mathbb{R}$, $t>0$, $\sigma$, $u$ and $\varepsilon$
 are stress, displacement and strain, as above,
 $\rho = \mathrm{const}.$ is the density of the media and $E=\mathrm{const}.$ is the Young
 modulus of elasticity. Equation (\ref{eq:4}) is the equilibrium equation and it
 is the consequence of the Second Newton Law. Equation (\ref{eq:5}) is the
 constitutive equation and in the case of elastic media it is known as the
 Hooke law. Equation (\ref{eq:6}) is the strain measure for the local small
 deformations, i.e. connection between strain and displacement. Wave equation
 is now obtained by substituting equation (\ref{eq:6}) into (\ref{eq:5}) and
 subsequently (\ref{eq:5}) into (\ref{eq:4}). It reads%
 \begin{equation} \label{eq:7}
 \frac{\partial^{2}}{\partial t^{2}}u(x,t) =
 c^{2}\frac{\partial^{2}}{\partial x^{2}}u(x,t), \qquad
 c=\sqrt{\frac{E}{\rho}}.
 \end{equation}
 Newly introduced parameter $c$ can be physically interpreted as the wave
 speed.

 The first authors who introduced the fractional Zener model were
 M. Caputo and F. Mainardi
 \cite{CaputoMainardi-1971a, CaputoMainardi-1971b}.
 Such model was also considered in
 \cite{Rabotnov} but without using fractional calculus.
 A similar problem of stress waves in a viscoelastic medium was
 investigated in \cite{GonsovskiRossikhin} by using the inversion
 formula of the classical Laplace transform.
 Wave equation (\ref{eq:7}) can be generalized
 in the framework of fractional calculus
 in several ways.
 One of the possibilities is to substitute
 second order partial derivative with respect to time by fractional one,
 as it has been done
 in \cite{Mainardi} for one dimension, and in \cite{Hanyga-stf} and
 \cite{Hanyga-tf} for $d$ dimensions, $d\in \{1,2,3\}$. Time partial derivative
 can also be replaced with either two fractional derivatives of different order,
 as in \cite{APZ07}, or by distributed order fractional derivative (integral over
 given range of order of fractional derivatives), as in
 \cite{Kochubei08} and \cite{MainardiPagniniGorenflo} for $\alpha \in (0,1)$,
 and in \cite{APZ09-I} and \cite{APZ09-II}
 for $\alpha \in [0,2]$. In each of these papers the intention
 was to stress similarity between wave and diffusion equation, since by
 allowing the order of fractional derivative to be $\alpha \in [0,2]$
 both equations can be obtained as special cases.

 In this article wave equation (\ref{eq:7}) is generalized for the case of
 viscoelastic media described by fractional Zener model. Generalization is
 conducted so that constitutive equation (\ref{eq:5}) is changed, since
 Hooke's
 law is inappropriate to describe viscoelastic media. Other two equations,
 namely equation of motion of deformable body (\ref{eq:4}) and relation between
 strain and displacement (\ref{eq:6}), are valid for any type of deformable
 body, with mentioned restrictions, therefore they remain unchanged. Similar
 generalization was done in \cite{RossikhinShitikova} in case of bounded domain.

 The paper is organized as follows. In Section \ref{sec:math prelim} we recall
 basic notions and fix notation which will be used throughout this paper. In
 Section \ref{sec:Cauchy problem} we introduce dimensionless quantities, derive
 (\ref{eq:system})
 and (\ref{eq:fzwe}) and set up the Cauchy problem (\ref{eq:fzwe}), (\ref{eq:ic})\ in
 the setting of distribution theory. Section \ref{sec:fund sol} is devoted to
 examination of the Cauchy problem (\ref{eq:fzwe}), (\ref{eq:ic}) where we prove a
 theorem on existence and uniqueness of its fundamental solution. In Section
 \ref{sec:num ex} we give a numerical example, which illustrates the behavior of
 the fundamental solution for the specific choice of the parameters $\alpha$
 and $\tau$, as well as initial conditions.

 \section{Mathematical preliminaries} \label{sec:math prelim}

 Let $\Omega$ be an open set in $\mathbb{R}^{n}$. The space of
 distributions will be denoted by $\cD'(\Omega)$ and the space
 of Schwartz tempered distributions by $\cS'(\mathbb{R}^{n})$.
 $\cD'_{+}(\mathbb{R})$ and $\cS'_{+}(\mathbb{R})$ are the subspaces of
 $\cD'(\mathbb{R})$ and $\cS'(\mathbb{R})$
 respectively, containing distributions supported on $[0,\infty)$.
 Since we shall mostly work with functions and distributions depending on
 two variables, $u=u(x,t)$, we introduce $\cS'(\mathbb{R}\times \mathbb{R}_{+})$
 to be the space of all distributions $u\in \cS'(\mathbb{R}^{2})$,
 which vanish for $t<0$.

 We shall also need the set $\mathbb{C}_{+}=\{z\in\mathbb{C}\,|\, \rep z>0\}$.

 Let $t\in [0,a]$, $a>0$, and $y\in AC([0,a])$. The left
 Riemann-Liouville fractional derivative of order $\alpha\in [0,1)$,
 ${}_{0}D_{t}^{\alpha}y$, is defined as
 $$
 {}_{0}D_{t}^{\alpha}y(t) =\frac{1}{\Gamma (1-\alpha)}
 \frac{d}{dt} \int_{0}^{t} \frac{y(\zeta)}{(t-\zeta)^{\alpha}}\,d\zeta,
 $$
 where $\Gamma$ is the Euler gamma function.

 In the distributional setting, one introduces a family
 $\{f_{\alpha}\}_{\alpha\in\mathbb{R}}\in \cD'_{+}(\mathbb{R})$
 as
 $$
 f_{\alpha}(t) = \left\{ \begin{array}{ll}
 H(t) \frac{t^{\alpha -1}}{\Gamma (\alpha)}, &
 \alpha >0, \\
 \frac{d^{N}}{dt^{N}}f_{\alpha +N}(t), & \alpha \leq 0, \alpha +N>0,
 N\in \mathbb{N}
 \end{array}
 \right. ,
 $$
 where $H$ is the Heaviside function. Then $f_{\alpha}\ast$ is a
 convolution operator acting on $\cD'_{+}(\mathbb{R})$
 (also, $f_{\alpha}\ast :\cS'_{+}(\mathbb{R})\to
 \cS'_{+}(\mathbb{R})$). For $\alpha <0$ it is called the
 operator of fractional differentiation. Moreover, for $y\in AC([0,a])$ it
 coincides with the left Riemann-Liouville fractional derivative, i.e.
 $$
 {}_{0}D_{t}^{\alpha}y=f_{-\alpha}\ast y.
 $$
 For $y\in \cS'(\mathbb{R})$ the Fourier transform is
 defined as
 $$
 \left\langle \cF y,\varphi \right\rangle =\left\langle y,\cF
 \varphi \right\rangle, \qquad \varphi \in \cS(\mathbb{R}),
 $$
 where for $\varphi \in \cS(\mathbb{R})$
 $$
 \cF\varphi (\xi)=\hat{\varphi}(\xi)=\int_{-\infty }^{\infty}
 \varphi (x)e^{-i\xi x}\,dx, \qquad \xi \in \mathbb{R}.
 $$
 Let $y\in \cD'_{+}(\mathbb{R})$ such that $e^{-\xi t}y\in
 \cS'(\mathbb{R})$, for all $\xi >a>0$. Then the Laplace
 transform of $y$ is defined by
 $$
 \cL y(s) =\widetilde{y}(s) = \cF(e^{-\xi t}y)(\eta),
 \qquad s=\xi +i\eta.
 $$
 It is well known that the function $\cL y$ is holomorphic in the half
 plane $\rep s>a$ (see e.g. \cite{DautryLions-vol5} or \cite{Vladimirov-EqMP}). In
 particular, for $y\in L^{1}(\mathbb{R})$ such that $y(t)=0$,
 for $t<0$, and $|y(t)| \leq Ae^{at}$ ($a,A>0$) the
 Laplace transform is
 $$
 \cL y(s)=\int_{0}^{\infty} y(t)e^{-st}\, dt,
 \qquad \rep s>0.
 $$
 We shall need the following properties: $y,y_{1},y_{2}\in
 \cS'(\mathbb{R})$
 \begin{eqnarray*}
 \cF\left[\frac{d^{n}}{dx^{n}}y\right] (\xi) &=& (i\xi)^{n}\cF y(\xi), \\
 \cL[{}_{0}D_{t}^{\alpha}y](s) &=& s^{\alpha}\cL y(s), \\
 \cL[y_{1}\ast y_{2}](s) & = & \cL y_{1}(s)\cdot \cL y_{1}(s), \\
 \cL\delta(s) &=& 1.
 \end{eqnarray*}%
 In order to introduce the inverse Laplace transform we recall from
 \cite{Vladimirov-EqMP} the following: Let $Y$ be a holomorphic function in the half
 plane $\rep s>a$ such that $|Y(s)| \leq A\frac{(1+|s|^{m})}{|\rep s|^k}$,
 $m,k\in \R$.
 Then there exists a distribution $y\in \cS'_{+}(\mathbb{R})$ such that
 $\cL y=Y$, and
 $$
 y(t)=\cL^{-1}Y(t) =\frac{1}{2\pi i}\int_{a-i\infty}^{a+i\infty}
 Y(s) e^{st}\, ds, \qquad t>0.
 $$

 The notion of fundamental solution is introduced as follows
 (see e.g. \cite{Vladimirov-EqMP}). Let $P$ be a
 linear partial integro-differential operator with constant coefficients. A
 fundamental solution of $P$, denoted by $E$, is a solution to the equation
 $Pu=\delta$. Once the fundamental solution is determined one finds a
 solution to $Pu=f$ as $u=E\ast f$.

 The Cauchy problem for the second order linear partial integro-differential
 operator with constant coefficients $P$ is given by
 \begin{gather}
 Pu(x,t)=f(x,t)  \label{eq:mp-cp} \\
 u(x,0)=u_{0}(x), \qquad
 \frac{\partial }{\partial t}u(x,0)=u_{1}(x),  \label{eq:mp-cp-ic}
 \end{gather}
 where $f$ is continuous for $t\geq 0$, $u_{0}\in C^{1}(\mathbb{R})$
 and $u_{1}\in C(\mathbb{R})$. A classical solution $u(x,t)$ to the
 Cauchy problem (\ref{eq:mp-cp}-\ref{eq:mp-cp-ic}) is of class
 $C^{2}$ for $t>0$, of class $C^{1}$ for $t\geq 0$, satisfies equation
 (\ref{eq:mp-cp}) for $t>0$, and initial conditions (\ref{eq:mp-cp-ic})
 when $t\to 0$. If functions $u$ and $f$ are continued by zero for
 $t<0$, then the following
 equation is satisfied in $\cD'(\mathbb{R}^2)$:
 \begin{equation} \label{eq:mp-gencp}
 Pu=f(x,t) +u_{0}(x) \delta'(t)
 +u_{1}(x) \delta (t).
 \end{equation}

 The classical solutions of the Cauchy problem
 (\ref{eq:mp-cp}-\ref{eq:mp-cp-ic}) are among
 those solutions of equation (\ref{eq:mp-gencp}) that vanish for $t<0$. Therefore,
 the problem of finding generalized solutions (in $\cD'(\mathbb{R}^2)$)
 of the equation (\ref{eq:mp-gencp}) that
 vanish for $t<0$ will be called {\it generalized Cauchy problem for the
 operator} $P$. If there is a fundamental solution $E$ of the operator $P$
 and if $f\in \cD'(\mathbb{R}^2)$
 vanishes for $t<0$ then there exists a unique solution to corresponding
 generalized Cauchy problem and is given by
 $$
 u=E\ast (f(x,t) +u_{0}(x) \delta'(t) +u_{1}(x) \delta(t)).
 $$

 Again, we refer to \cite{DautryLions-vol5}, \cite{Treves-PDEs} or
 \cite{Vladimirov-EqMP} for more details.

 \section{Set up of the Cauchy problem} \label{sec:Cauchy problem}

 In order to obtain (\ref{eq:system}), we introduce dimensionless coordinates in the
 system (\ref{eq:4}-\ref{eq:6}), where (\ref{eq:5}) is modified by the fractional
 Zener model of the viscoelastic body:
 \begin{equation} \label{eq:frac Zener model}
 \sigma (x,t) +\tau_{\sigma} {}_{0}D_{t}^{\alpha}\sigma (x,t)
 =E[\varepsilon (x,t) +\tau_{\varepsilon} {}_{0}D_{t}^{\alpha}\varepsilon (x,t)],
 \qquad x\in \mathbf{\mathbb{R}}, t>0,
 \end{equation}
 where $\tau_{\sigma}$, $\tau_{\varepsilon}$ are relaxation times
 satisfying $\tau_{\varepsilon}>\tau_{\sigma}>0$. The latter condition
 follows from the Second Law of Thermodynamics (see \cite{Atanackovic02}).

 We subject to the system (\ref{eq:4}), (\ref{eq:6}) and
 (\ref{eq:frac Zener model}) initial conditions
 \begin{gather*}
 u(x,0) =u_{0}(x), \qquad \frac{\partial}{\partial t} u(x,0) =v_{0}(x), \\
 \sigma(x,0) =0, \qquad \varepsilon (x,0) =0,
 \end{gather*}
 where $x\in \mathbf{\mathbb{R}}$, $u_{0}$ and $v_{0}$ are initial
 displacement and velocity. Note that
 there are no initial stress and strain. We also supply boundary conditions
 $$
 \lim_{x\to \pm \infty}u(x,t)=0, \qquad \lim_{x\to \pm \infty}\sigma (x,t)=0,
 \qquad t\geq 0.
 $$

 We will need the following lemma.

 \blem
 Let $y\in AC([0,a])$, $a>0$, $0\leq \alpha <1$ and
 $T^{\ast}>0$. Let $t$ and $y$ be transformed as $(\bar{t},\bar{y})
 =\left(\frac{t}{T^{\ast}},y\right)$, then the left
 Riemann-Liouville fractional derivative ${}_{0}D_{t}^{\alpha}y(t)$,
 $t\geq 0$, is transformed as follows
 $$
 {}_{0}D_{\bar{t}}^{\alpha}\bar{y}(\bar{t})=(T^{\ast})^{\alpha}
 {}_{0}D_{t}^{\alpha}y(t).
 $$
 \elem

 \pr
 Since $\bar{y}(\bar{t})=y(t) =y(\bar{t}T^{\ast})$,
 we have
 \begin{eqnarray*}
 {}_{0}D_{\bar{t}}^{\alpha}\bar{y}(\bar{t}) &=&
 \frac{1}{\Gamma(1-\alpha)} \frac{d}{d\bar{t}}
 \int_{0}^{\bar{t}} \frac{\bar{y}(\xi)}{(\bar{t}-\xi)^{\alpha}}\,d\xi \\
 &=& \frac{1}{\Gamma(1-\alpha )}T^{\ast}\frac{d}{dt}
 \int_{0}^{\frac{t} {T^{\ast}}}\frac{y(\xi T^{\ast})}{(\frac{t}{T^{\ast}}-\xi)^{\alpha}}\,d\xi \\
 &=& (T^{\ast})^{\alpha} \frac{1}{\Gamma(1-\alpha)} \frac{d}{dt}
 \int_{0}^{t} \frac{y(z)}{(t-z)^{\alpha}}\,dz \\
 &=& (T^{\ast})^{\alpha}{}_{0}D_{t}^{\alpha}y(t).
 \end{eqnarray*}
 \ep

 Let $0\leq \alpha <1$, $\tau_{\varepsilon}>\tau_{\sigma }>0$ and $E$ be
 constants appearing in (\ref{eq:frac Zener model}),
 and let $\rho$ be the constant in (\ref{eq:4}).
 Set $T^{\ast}=\sqrt[\alpha]{\tau_{\varepsilon}}$ and
 $X^{\ast }=\sqrt[\alpha]{\tau_{\varepsilon}}\sqrt{\frac{E}{\rho}}$.
 Then $X^{\ast}$ and $T^{\ast}$ are constants (measured in meters and seconds
 respectively) which
 satisfy $\left(\frac{X^{\ast}}{T^{\ast}}\right)^{2}\frac{\rho}{E}=1$
 and $(T^{\ast})^{\alpha}=\frac{1}{\tau_{\varepsilon}}$.

 We introduce dimensionless quantities in (\ref{eq:4}),
 (\ref{eq:frac Zener model}) and (\ref{eq:6}) as
 $$
 \bar{x}=\frac{x}{X^{\ast}};
 \quad
 \bar{t}=\frac{t}{T^{\ast}};
 \quad
 \bar{u}=\frac{u}{X^{\ast}};
 \quad
 \bar{\sigma}=\frac{\sigma}{E}
 $$
 and
 $$
 \bar{u}_{0}=\frac{u_{0}}{X^{\ast}};
 \quad
 \bar{v}_{0}=\frac{v_{0}}{X^{\ast}}T^{\ast};
 \quad
 \bar{\tau}=\frac{\tau_{\sigma}}{\tau_{\varepsilon}}
 $$
 and obtain system (\ref{eq:system}-\ref{eq:bc}). Note that
 $\varepsilon$ is already dimensionless quantity.

 In order to simplify notation, a bar $\bar{\quad}$ will be dropped and
 dimensionless quantities: $\bar{x}$, $\bar{t}$, $\bar{\sigma}$, $\bar{u},$
 $\bar{\tau}$, $\bar{u}_{0}$ and $\bar{v}_{0}$ will be written as: $x, t,
 \sigma, u, \tau, u_{0}$ and $v_{0}$. Obviously, condition
 $0<\tau_{\sigma}<\tau_{\varepsilon}$ implies $0<\tau <1$.

 In our work an appropriate setting for studying system (\ref{eq:system}-\ref{eq:bc})
 will be distributional one. In fact, we will look for a fundamental solution
 to the generalized Cauchy problem for FZWE (\ref{eq:fzwe}) in
 $\cS'(\mathbb{R}\times \mathbb{R}_{+})$. This suffices to obtain a solution to
 (\ref{eq:system}-\ref{eq:bc}), since (\ref{eq:system}) and (\ref{eq:fzwe}) are equivalent.
 Indeed, by applying the Laplace transform with respect to time variable $t$
 to the second equation in (\ref{eq:system}), one obtains
 $$
 (1+\tau s^{\alpha})\widetilde{\sigma}(x,s)=
 (1+s^{\alpha})\widetilde{\varepsilon}(x,s).
 $$
 According to \cite{Oparnica02},
 $\cL^{-1}\left(\frac{1+s^{\alpha}}{1+\tau s^{\alpha}}\right)$
 is well-defined element in $\cS'_{+}(\mathbb{R})$, hence
 \begin{equation} \label{eq:sol sigma}
 \sigma =\cL^{-1}\left(\frac{1+s^{\alpha}}{1+\tau s^{\alpha}}\right)
 \ast_{t}\varepsilon.
 \end{equation}
 Inserting $\varepsilon$ from the third equation in (\ref{eq:system}) into
 (\ref{eq:sol sigma}) and then $\sigma$ into the first equation in
 (\ref{eq:system}) yields
 \begin{equation} \label{eq:fzwe2}
 \frac{\partial^{2}}{\partial t^{2}}u(x,t) =\cL^{-1}
 \left(\frac{1+s^{\alpha}}{1+\tau s^{\alpha}}\right)
 \ast_{t}
 \frac{\partial^{2}}{\partial x^{2}}u(x,t).
 \end{equation}
 Setting $L(t)=\cL^{-1}\left(\frac{1+s^{\alpha}}{1+\tau s^{\alpha}}\right)
 \ast_{t}$ we come to (\ref{eq:fzwe}). Therefore we have proved that
 (\ref{eq:system}) and (\ref{eq:fzwe}) are equivalent. Notice that equation
 (\ref{eq:fzwe2}) is of the form $Pu=0$, with
 \begin{equation} \label{eq:op P}
 P:=\frac{\partial^{2}}{\partial t^{2}}-\cL^{-1}
 \left(\frac{1+s^{\alpha}}{1+\tau s^{\alpha}}\right)
 \ast_{t}
 \frac{\partial^{2}}{\partial x^{2}}.
 \end{equation}

 \section{The existence and uniqueness of a solution to the Cauchy problem
 (\ref{eq:fzwe}), (\ref{eq:ic})} \label{sec:fund sol}

 The aim of this section is to find a solution to the generalized Cauchy
 problem to (\ref{eq:fzwe}), i.e.
 \begin{equation} \label{eq:fzwe-distr}
 \frac{\partial^{2}}{\partial t^{2}}u(x,t) =
 \cL^{-1}\left(\frac{1+s^{\alpha}}{1+\tau s^{\alpha}}\right)
 \ast_{t} \frac{\partial^{2}}{\partial x^{2}}u(x,t)
 + u_{0}(x) \delta'(t)+v_{0}(x)\delta (t),
 \end{equation}
 or equivalently
 $$
 Pu(x,t)=u_{0}(x)\delta'(t)+v_{0}(x)\delta (t),
 $$
 where $P$ is given by (\ref{eq:op P}).

 \brem
 We have already explained that initial conditions are included
 into the generalized Cauchy problem (see Section \ref{sec:math
 prelim}). For functions in $\cS'$ boundary conditions
 (\ref{eq:bc}) are automatically fulfilled.
 \erem

 We state the main theorem.

 \bthm \label{th:glavna}
 Let $u_{0},v_{0}\in \cS'(\mathbb{R})$. Then
 there exists a unique solution $u\in\cS'(\mathbb{R}\times \mathbb{R}_{+})$
 to (\ref{eq:fzwe-distr}) given by
 \begin{equation} \label{eq:sol-u}
 u(x,t)=S(x,t) \ast_{x,t}(u_{0}(x)\delta'(t)+v_{0}(x)\delta (t)),
 \end{equation}
 where
 \begin{eqnarray}
 S(x,t) &=& 1+\frac{1}{4\pi i}
 \int_{0}^{\infty}
 \left(
 \sqrt{\frac{1+\tau q^{\alpha}e^{i\alpha \pi}}{1+q^{\alpha}e^{i\alpha \pi}}}
 e^{|x| q
 \sqrt{\frac{1+\tau q^{\alpha}e^{i\alpha \pi}}{1+q^{\alpha}e^{i\alpha \pi}}}}\right.
 \nonumber \\
 && \qquad
 \left. -\sqrt{\frac{1+\tau q^{\alpha}e^{-i\alpha \pi}}{1+q^{\alpha}e^{-i\alpha \pi}}}
 e^{|x| q
 \sqrt{\frac{1+\tau q^{\alpha}e^{-i\alpha \pi}}{1+q^{\alpha}e^{-i\alpha \pi}}}}\right)
 \frac{e^{-qt}}{q}\, dq,  \label{eq:fund sol}
 \end{eqnarray}
 is the fundamental solution of operator $P$,
 $S\in\cS'(\mathbb{R}\times\mathbb{R}_{+})$ with support in the
 cone $|x| <\frac{t}{\sqrt{\tau}}$.
 \ethm

 We will need the following lemma.

 \blem \label{lem:ode}
 Let $f\in \cS'(\mathbb{R})$. Then the equation
 \begin{equation} \label{eq:we1}
 v''-\omega v=-f
 \end{equation}
 has a solution $v\in \cS'(\mathbb{R})$ for all $\omega \in
 \mathbb{C}\setminus (-\infty ,0]$, which is of the form
 $$
 v = \frac{e^{-\sqrt{\omega }|x|}}{2\sqrt{\omega}}\ast f,
 $$
 where $\sqrt{\omega}$ is the main branch.
 \elem

 \pr
 We first apply Fourier transform to (\ref{eq:we1}):
 \begin{equation} \label{eq:Fwe1}
 \widehat{v}=\frac{1}{\xi^{2}+\omega}\widehat{f}.
 \end{equation}
 In order to find $v$ we need to calculate
 $$
 \cF^{-1}\left(\frac{1}{\xi^{2}+\omega}\right) (x)
 = \frac{1}{2\pi}
 \int_{\mathbb{R}} e^{ix\xi}\frac{1}{\xi^{2}+\omega}\,d\xi.
 $$
 Let $x\geq 0$. Let $\Gamma_{+}=\Gamma_{1}\cup \Gamma_{2}$, where
 $\Gamma_{1}=\{z=Re^{i\varphi}\;|\;R>0,0<\varphi<\pi\}$ and
 $\Gamma_{2}=\{z=\xi\;|\;-R<\xi <R\}$.
 By the Cauchy residual theorem we have that
 $$
 \oint_{\Gamma} e^{ixz}\frac{1}{z^{2}+\omega}\,dz
 = 2\pi i \sum_{i=1}^{k} \mbox{Res}(z_{i}),
 $$
 where $z_{i}$ are poles of the function
 $z\mapsto \frac{e^{ixz}}{z^{2}+\omega}$
 lying inside of $\Gamma$. For $\omega \in \mathbb{C}\setminus (-\infty ,0]$
 there is only one such pole $z_{1}=i\sqrt{\omega}$ and
 $\mbox{Res}(z_{1})= \left. \frac{e^{ixz}}{2z}\right\vert_{z=z_{1}}
 =\frac{e^{-x\sqrt{\omega}}}{2i\sqrt{\omega}}$. Letting
 $R\to \infty$, the integral over $\Gamma_{1}$ tends to zero and
 $\int_{\Gamma_{2}} e^{ixz}\frac{1}{z^{2}+\omega}\,dz
 = \int_{\mathbb{R}} e^{ix\xi}\frac{1}{\xi^{2}+\omega}\,d\xi$, thus
 $$
 \cF^{-1}\left(\frac{1}{\xi^{2}+\omega}\right)(x) =
 \frac{1}{2\pi} \int_{\Gamma_{2}} e^{ixz}\frac{1}{z^{2}+\omega}\,dz =
 \frac{e^{-x\sqrt{\omega}}}{2\sqrt{\omega}}.
 $$
 For $x<0$ and $\Gamma_{-}=\Gamma_{3}\cup \Gamma_{eq:4}$, where
 $\Gamma_{3}=\{z=Re^{i\varphi}\;|\;R>0,-\pi <\varphi <0\}$ and
 $\Gamma_{4}=\{z=-\xi \;|\;-R<\xi <R\}$ we can apply the same arguments as for the
 case $x\geq 0$ to obtain $\cF^{-1}\left(\frac{1}{\xi^{2}+\omega}\right) (x)
 =\frac{e^{x\sqrt{\omega}}}{2\sqrt{\omega}}$.

 Thus for $x\in\mathbb{R}$ we can write
 $$
 \cF^{-1}\left(\frac{1}{\xi^{2}+\omega}\right)(x)
 = \frac{e^{-\sqrt{\omega}|x|}}{2\sqrt{\omega}}.
 $$
 Now, the inverse Fourier transform applied to (\ref{eq:Fwe1}) yields the result.
 \ep

 {\bf Proof of the main theorem.}
 Applying Laplace transform with respect to $t$
 to (\ref{eq:fzwe-distr}) we obtain:
 \begin{equation} \label{eq:lt-fzwe-distr}
 \frac{\partial^{2}}{\partial x^{2}}\widetilde{u}(x,s)
 - s^{2} \frac{1+\tau s^{\alpha}}{1+s^{\alpha}}\widetilde{u}(x,s) =
 -\frac{1+\tau s^{\alpha}}{1+s^{\alpha}}(su_{0}(x)+v_{0}(x)).
 \end{equation}
 Set $\omega (s):=s^{2}\frac{1+\tau s^{\alpha}}{1+s^{\alpha}}$ and
 $f(x,s):= \frac{1+\tau s^{\alpha}}{1+s^{\alpha}}(su_{0}(x)+v_{0}(x))$,
 $s\in \mathbb{C}_{+}$.
 In order to apply Lemma \ref{lem:ode} to (\ref{eq:lt-fzwe-distr}) we
 have to show
 first that $f(\cdot,s)\in \cS'(\mathbb{R})$, which follows from
 assumptions that $u_{0},v_{0}\in \cS'(\mathbb{R})$ and
 second, that
 $\omega (s)\in \mathbb{C}\setminus (-\infty ,0]$, for all
 $s\in \mathbb{C}_{+}$. The latter can be verified in the following way.
 Take arbitrary $s=\rho e^{i\varphi}$, $\rho>0$, $-\frac{\pi}{2}<\varphi
 <\frac{\pi}{2}$.
 Then $\omega (s)=\rho^{2}e^{2i\varphi}\frac{1+\tau \rho^{\alpha}
 e^{i\alpha \varphi}}{1+\rho^{\alpha}e^{i\alpha \varphi}}$, and after a
 straightforward calculation we arrive to
 \begin{equation} \label{eq:Re oms}
 \rep\omega (s)= \frac{\rho^{2}}{A}
 \left[\left(
 1+\rho^{\alpha}(1+\tau) \cos(\alpha \varphi) +\tau \rho^{2\alpha}
 \right)
 \cos(2\varphi) +\rho^{\alpha}(1-\tau) \sin(\alpha \varphi)\sin(2\varphi)
 \right],
 \end{equation}
 \begin{equation} \label{eq:Im oms}
 \imp\omega (s)=\frac{\rho^{2}}{A}
 \left[ \left(
 1+\rho^{\alpha}(1+\tau) \cos(\alpha \varphi) +\tau \rho^{2\alpha}
 \right)
 \sin(2\varphi) -\rho^{\alpha}(1-\tau)\sin(\alpha \varphi) \cos(2\varphi)
 \right],
 \end{equation}
 where $A=(1+\rho^{\alpha}\cos(\alpha \varphi))^{2}
 +\rho^{2\alpha}\sin^{2}(\alpha \varphi)>0$. Suppose
 that $\omega (s)\in (-\infty ,0]$, for some $s\in \mathbb{C}_{+}$. Then
 $\imp\omega (s) =0$, and (\ref{eq:Im oms}) yields
 \begin{equation} \label{eq:kontra}
 \left(1+\rho^{\alpha}(1+\tau) \cos(\alpha \varphi)
 +\tau \rho^{2\alpha}\right) =
 \rho^{\alpha}(1-\tau) \sin(\alpha \varphi)
 \frac{\cos(2\varphi)}{\sin(2\varphi)}.
 \end{equation}
 From (\ref{eq:Re oms}) we now obtain
 $$
 \rep\omega (s)=\frac{\rho^{2}}{A}\rho^{\alpha}(1-\tau)
 \sin(\alpha \varphi) \sin(2\varphi)
 \left( \frac{\cos^{2}(2\varphi)}{\sin^{2}(2\varphi)}+1 \right)
 \leq 0.
 $$
 where we have used (\ref{eq:kontra}) and the assumption $\omega(s)
 \in (-\infty ,0]$, for all $s\in \mathbb{C}_{+}$.

 Since $\frac{\rho^{2}}{A}\rho^{\alpha}(1-\tau) \sin(\alpha \varphi)
 \left( \frac{\cos^{2}(2\varphi)}{\sin^{2}(2\varphi)}+1\right) >0$,
 for $\rho >0$, $-\frac{\pi}{2}<\varphi <\frac{\pi }{2}$,
 it follows that $\sin(2\varphi) <0$,
 and therefore $\varphi \in \left(-\frac{\pi}{2},-\frac{\pi}{4}\right]
 \cup \left[\frac{\pi}{},\frac{\pi}{2}\right)$. However, (\ref{eq:kontra})
 can not be satisfied for those $\varphi$.
 (Indeed, for $\varphi \in \left(\frac{\pi}{2},-\frac{\pi}{4}\right]$,
 $\cos(\alpha \varphi), \cos(2\varphi), \sin(\alpha \varphi) >0$, but
 $\sin(2\varphi) <0$ and similarly for $\varphi\in\left[\frac{\pi}{4},
 \frac{\pi}{2}\right)$.) Hence, $\omega (s)\in \mathbb{C}\setminus
 (-\infty ,0]$, for all $s\in \mathbb{C}_{+}$.

 Now we can apply Lemma \ref{lem:ode} to obtain a solution of
 (\ref{eq:lt-fzwe-distr}):
 \begin{eqnarray}
 \widetilde{u}(x,s) &=&
 \frac{e^{-\sqrt{\omega(s)} |x|}}{2\sqrt{\omega(s)}}
 \ast_{x}
 \frac{\omega (s)}{s^{2}}(su_{0}(x)+v_{0}(x)) \nonumber \\
 &=& \frac{\sqrt{\omega (s)}e^{-\sqrt{\omega(s)} |x|}}{2s^{2}}
 \ast_{x}
 (su_{0}(x)+v_{0}(x)).  \label{eq:tilde u}
 \end{eqnarray}%
 To obtain a solution $u$ to (\ref{eq:fzwe}) we need to calculate inverse
 Laplace transform of (\ref{eq:tilde u}). For that purpose set
 \begin{equation} \label{eq:tilde S}
 \widetilde{S}(x,s):=\frac{\sqrt{\omega (s)}e^{-\sqrt{\omega (s)}|x|}}{2s^{2}}
 =\frac{1}{2s}\sqrt{\frac{1+\tau s^{\alpha}}{1+s^{\alpha}}}
 e^{-|x|s\sqrt{\frac{1+\tau s^{\alpha}}{1+s^{\alpha}}}}, \qquad
 x\in \mathbb{R}, s\in \mathbb{C}_{+},
 \end{equation}
 and
 $$
 S(x,t) =\cL^{-1}\left[\widetilde{S}(x,s)\right](t), \qquad
 x\in\mathbb{R}, t>0.
 $$
 $S(x,t)$ is well-defined since
 $|\widetilde{S}(x,s)|\leq C\frac{(1+|s|^m)}{|\rep s|}$, $m\in\mathbb{N}$,
 so its inverse Laplace transform exists (cf.\ Section
 \ref{sec:math prelim}).

 Note also that $S$ is a fundamental solution of $P$.
 Then the solution $u$ is given by (\ref{eq:sol-u}).

\begin{figure}[htb]
 \begin{center}
 \includegraphics[width=6cm]{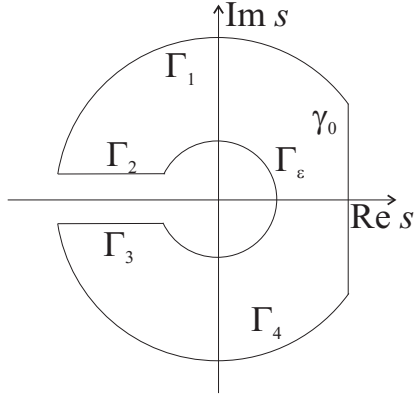}
 \end{center}
 \caption{Integration contour $\Gamma$}
 \label{fig-1}
 \end{figure}

 Multiform function $\widetilde{S}$, given by (\ref{eq:tilde S}), has
 branch points at $s=0$ and $s=\infty$ and has no singularities.
 Hence, it can be evaluated by the use of the Cauchy integral formula
 \begin{equation} \label{eq:Cauchy int for}
 \oint_{\Gamma} \widetilde{S}(x,s) e^{st}\,ds=0, \qquad
 x\in \mathbb{R}, t>0,
 \end{equation}
 where $\Gamma =\Gamma_{1}\cup \Gamma_{2}\cup \Gamma_{\varepsilon}\cup
 \Gamma_{3}\cup \Gamma_{4}\cup\gamma_{0}$, is a contour given in Figure
 \ref{fig-1}.

 More precisely, for arbitrarily chosen $R>0$, $0<\varepsilon <R$ and $a>0$,
 $\Gamma$ is defined by
 $$
 \begin{array}{rl}
 \Gamma_{1}: & s=Re^{i\varphi},\varphi_{0}<\varphi <\pi; \\
 \Gamma_{2}: & s=qe^{i\pi},-R<-q<-\varepsilon; \\
 \Gamma_{\varepsilon}: & s=\varepsilon e^{i\varphi},-\pi <\varphi <\pi; \\
 \Gamma_{3}: & s=qe^{-i\pi},\varepsilon <q<R; \\
 \Gamma_{4}: & s=Re^{i\varphi},-\pi <\varphi<-\varphi_{0}; \\
 \gamma_{0}: & s=a(1+i\tan \varphi),-\varphi_{0}<\varphi <\varphi_{0},
 \end{array}
 $$
 where $\varphi_{0}=\arccos(\frac{a}{R})$. Note that $\lim_{R\to\infty}
 \varphi_{0}=\frac{\pi}{2}$. In the limit when $R\to \infty$,
 integral along contour $\Gamma_{1}$ reads ($x\in \mathbb{R}, t>0$)
 \begin{equation} \label{eq:int Gamma1}
 \lim_{R\to\infty} \int_{\Gamma_{1}}
 \widetilde{S}(x,s) e^{st}\,ds =
 \frac{1}{2} \lim_{R\to \infty} \int_{\varphi_{0}}^{\pi}
 \sqrt{\frac{1+\tau R^{\alpha}e^{i\alpha\varphi}}{1+R^{\alpha}e^{i\alpha \varphi}}}
 e^{-|x|Re^{i\varphi}\sqrt{\frac{1+\tau R^{\alpha}e^{i\alpha\varphi}}{1+R^{\alpha}e^{i\alpha \varphi}}}
 +Rte^{i\varphi}}i\,d\varphi.
 \end{equation}
 Evaluating the absolute value of $\int_{\Gamma_{1}}\widetilde{S}(x,s)
 e^{st}\,ds$ one obtains ($x\in \mathbb{R}, t>0$)
 $$
 \lim_{R\to \infty} \left\vert \int_{\Gamma_{1}}
 \widetilde{S}(x,s) e^{st}\,ds \right\vert \leq
 \frac{1}{2} \lim_{R\to\infty} \int_{\varphi_{0}}^{\pi}
 \left\vert
 \sqrt{\frac{1+\tau R^{\alpha}e^{i\alpha \varphi}}{1+R^{\alpha}e^{i\alpha \varphi}}}
 \right\vert
 \left\vert e^{-|x| Re^{i\varphi}\sqrt{\frac{1+\tau R^{\alpha}e^{i\alpha\varphi}}{1+R^{\alpha}e^{i\alpha\varphi}}}}
 \right\vert
 e^{Rt\cos\varphi}\,d\varphi.
 $$
 In the limit when $R\to \infty$ the expression
 $\sqrt{\frac{1+\tau R^{\alpha}e^{i\alpha \varphi}}{1+R^{\alpha}e^{i\alpha \varphi}}}$
 tends to $\sqrt{\tau}$ and therefore
 $$
 \lim_{R\to \infty} \left\vert \int_{\Gamma_{1}}
 \widetilde{S}(x,s) e^{st}\,ds
 \right\vert
 \leq
 \frac{\sqrt{\tau}}{2}
 \lim_{R\to \infty} \int_{\varphi_{0}}^{\pi}
 e^{R\cos \varphi (t-|x| \sqrt{\tau})}\,d\varphi =0,
 \qquad \mbox{ if }
 t>|x| \sqrt{\tau},
 $$
 since $\cos \varphi <0$ for $\varphi \in
 \left(\frac{\pi}{2},\pi\right)$.
 Similar argument is valid for the integral along $\Gamma_{4}$.

 In the limit when $\varepsilon \to 0$, the integral along
 $\Gamma_{\varepsilon}$ is given by the formula similar to
 (\ref{eq:int Gamma1}) and it is calculated as
 $$
 \lim_{\varepsilon \to 0} \int_{\Gamma_{\varepsilon}}
 \widetilde{S}(x,s) e^{st}\, ds =
 \frac{1}{2} \lim_{\varepsilon \to 0} \int_{\pi}^{-\pi}
 \sqrt{\frac{1+\tau \varepsilon^{\alpha}e^{i\alpha \varphi}}{1+\varepsilon^{\alpha}e^{i\alpha\varphi}}}
 e^{-|x| \varepsilon e^{i\varphi}\sqrt{\frac{1+\tau\varepsilon^{\alpha}e^{i\alpha\varphi}}{1+\varepsilon^{\alpha}e^{i\alpha\varphi}}}
 +\varepsilon te^{i\varphi}}i\,d\varphi = -i\pi.
 $$
 Integrals along contours: $\Gamma_{2}$, $\Gamma_{3}$ and $\gamma_{0}$, in
 the limit when $R\to\infty$, $\varepsilon \rightarrow 0$, give
 \begin{eqnarray*}
 \lim_{\substack{R\to\infty \\ \varepsilon \to 0}}
 \int_{\Gamma_{2}}
 \widetilde{S}(x,s) e^{st}\, ds
 &=&
 -\frac{1}{2} \int_{0}^{\infty}
 \sqrt{\frac{1+\tau q^{\alpha}e^{i\alpha \pi}}{1+q^{\alpha}e^{i\alpha \pi}}}
 e^{-q\left(t-|x|\sqrt{\frac{1+\tau q^{\alpha}e^{i\alpha \pi}}{1+q^{\alpha}e^{i\alpha \pi}}}\right)}
 \frac{dq}{q} \\
 \lim_{\substack{R\to \infty \\ \varepsilon \to 0}}
 \int_{\Gamma_{3}}
 \widetilde{S}(x,s) e^{st}\,ds
 &=&
 \frac{1}{2} \int_{0}^{\infty}
 \sqrt{\frac{1+\tau q^{\alpha}e^{-i\alpha \pi}}{1+q^{\alpha}e^{-i\alpha \pi}}}
 e^{-q\left(t-|x| \sqrt{\frac{1+\tau q^{\alpha}e^{-i\alpha \pi}}{1+q^{\alpha}e^{-i\alpha \pi}}}\right)}
 \frac{dq}{q} \\
 \lim_{R\to \infty} \int_{\gamma_{0}}
 \widetilde{S}(x,s) e^{st}\,ds
 &=&
 2\pi i S(x,t).
 \end{eqnarray*}
 Now, by the Cauchy integral formula (\ref{eq:Cauchy int for}), we obtain $S$
 as in (\ref{eq:fund sol}). Therefore we have proved the theorem.
 \ep

 As a consequence of what we proved in Sections \ref{sec:Cauchy problem}
 and \ref{sec:fund sol}, we have the following corollary.

 \bcor
 Let $u$ be given by (\ref{eq:sol-u}). Then
 $$
 (u,\varepsilon,\sigma)(x,t)=\left( u(x,t),\frac{\partial}{\partial x}
 u(x,t),\cL^{-1}\left(\frac{1+s^{\alpha}}{1+\tau s^{\alpha}}\right)
 \ast_{t}\frac{\partial}{\partial x}u(x,t)\right) \in
 (\cS'(\mathbb{R}\times \mathbb{R}_{+}))^{3},
 $$
 is a unique solution to the system (\ref{eq:system}-\ref{eq:bc}).
 \ecor

 In a similar way one can also consider nonhomogeneous case.

 \brem
 Let the first equation in (\ref{eq:system}) be replaced by
 $$
 \frac{\partial}{\partial x}\sigma (x,t) =
 \frac{\partial^{2}}{\partial t^{2}} u(x,t) + f(x,t),
 $$
 where $f\in \cS'(\mathbb{R}\times \mathbb{R}_{+})$. This is a case
 of a rod under the influence of body forces. Then, the solution of the
 generalized Cauchy problem
 $$
 \frac{\partial^{2}}{\partial t^{2}}u(x,t) =
 \cL^{-1}\left(\frac{1+s^{\alpha}}{1+\tau s^{\alpha}}\right)
 \ast_{t}
 \frac{\partial^{2}}{\partial x^{2}}u(x,t) + f(x,t)+u_{0}(x)\delta'(t)
 +v_{0}(x)\delta (t),
 $$
 is
 $$
 u(x,t)=S(x,t) \ast_{x,t} (f(x,t)+u_{0}(x)\delta'(t)+v_{0}(x)\delta (t))
 $$
 where $S$ is as in (\ref{eq:fund sol}).
 \erem

 \brem
 Dimensionless condition $|x| <\frac{t}{\sqrt{\tau}}$,
 or $|x| <\sqrt{\frac{\rho}{E}}\sqrt{\frac{\tau_{\sigma}}{\tau_{\varepsilon}}}t$
 in the dimensional form, can be physically
 interpreted as the wave property. Namely, in the moment $t$, wave caused by
 the initial disturbance $u_{0}(x) =\delta (x)$, $x\in \mathbb{R}$,
 has reached the point which is at distance $|x|$
 from the origin, where the initial disturbance was applied. Thus, constant
 $c_{Z}=\sqrt{\frac{\rho}{E}}\sqrt{\frac{\tau_{\sigma}}{\tau_{\varepsilon}}}$
 can be interpreted as the wave speed in the fractional viscoelastic
 media of Zener type.
 \erem

 \section{A numerical example} \label{sec:num ex}

\begin{figure}[htb]
 \centerline{\includegraphics[width=6cm]{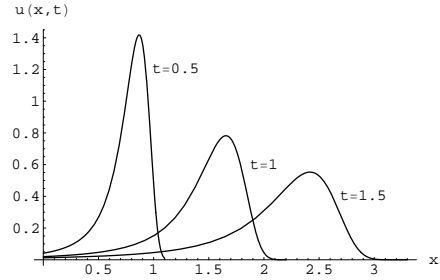}}
 \caption{Solution $u(x,t)$, $x\in (0,3)$, $t\in \{0.5,1,1.5\}$}
 \label{fig-2}
 \end{figure}

 Let $u_{0}=\delta$ and $v_{0}=0$ in the solution (\ref{eq:sol-u}) to the Cauchy
 problem of wave equation for fractional Zener type viscoelastic media. Then
 solution reads
 \begin{equation} \label{eq:num sol}
 u(x,t) = \frac{\partial}{\partial t}S(x,t),
 \qquad x\in\mathbb{R}, t>0,
 \end{equation}
 i.e. the fundamental solution represents the solution itself. Figure
 \ref{fig-2} presents the plots of $u(x,t)$,
 $x\in (0,3)$, $t\in \{0.5,1,1.5\}$, given by (\ref{eq:num sol}), for the set of
 parameters: $\alpha =0.23$ and $\tau =0.004$.

 In order to see the effect of a change in the order of a
 fractional derivative, we present following figures. Figure
 \ref{fig-3} presents plots of $u(x,t)$, $x\in (0,3)$, $t\in
 \{0.5,1,1.5\}$, given by (\ref{eq:num sol}), for different values
 of parameter $\alpha$. Plot of $u$, denoted by the dashed line
 corresponds to $\alpha =0.25$, dot-dashed line is used for $\alpha
 =0.5$, while the solid line denotes the plot for $\alpha =0.75$.
 Parameter $\tau =0.004$ is the same in all three figures.

 \begin{figure}[htb]
 \begin{center}$
 \begin{array}{cc}
 \includegraphics[scale=0.5]{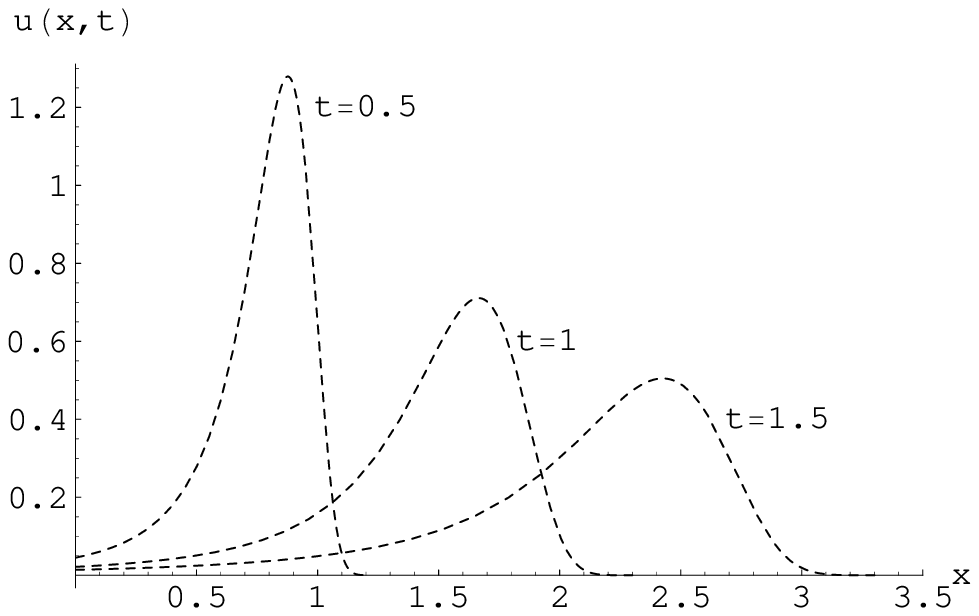}
 \includegraphics[scale=0.5]{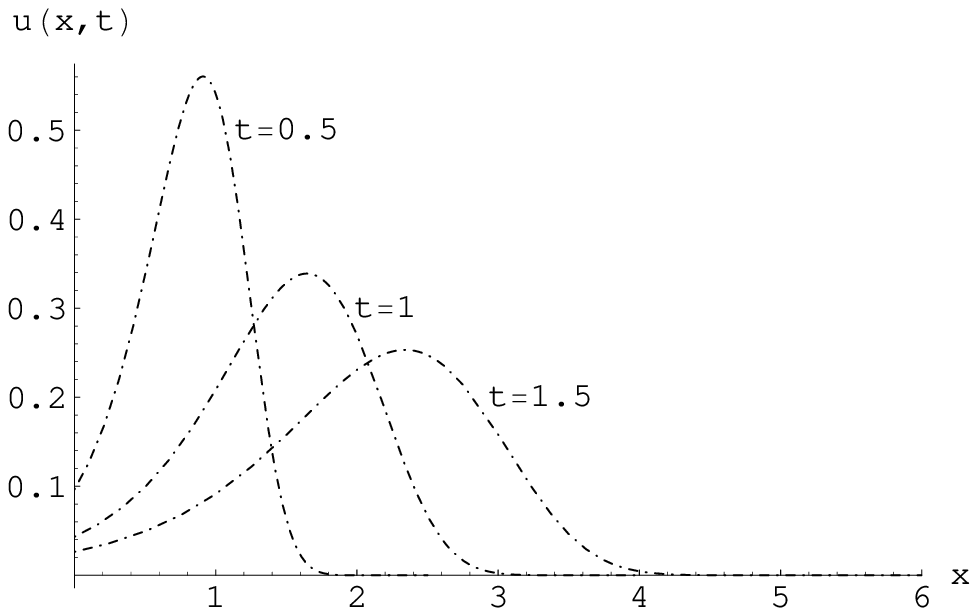}
 \includegraphics[scale=0.5]{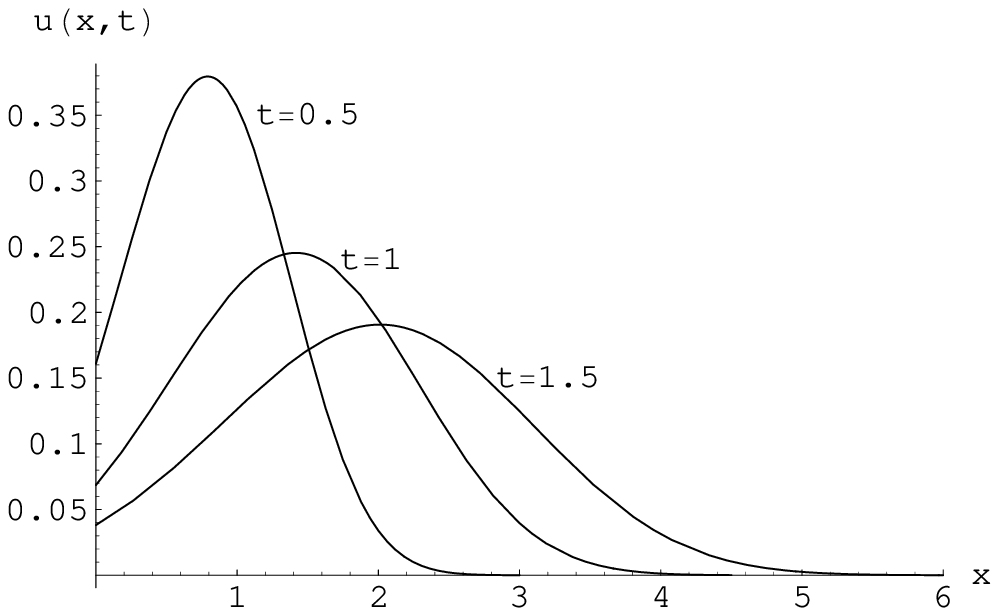}
 \end{array}$
 \end{center}
 \caption{Solution $u(x,t)$, $x\in (0,6)$, $t\in \{0.5,1,1.5\}$,
 $\alpha \in \{0.25,0.5,0.75\}$} \label{fig-3}
 \end{figure}

 From figure \ref{fig-3} one can see that for each value of the order of
 fractional derivative $\alpha$, as time increases, the maximum
 value of $u$ decreases, which is the consequence of the
 dissipative model of a media.

\begin{figure}[htb]
 \begin{center}$
 \begin{array}{cc}
 \includegraphics[scale=0.5]{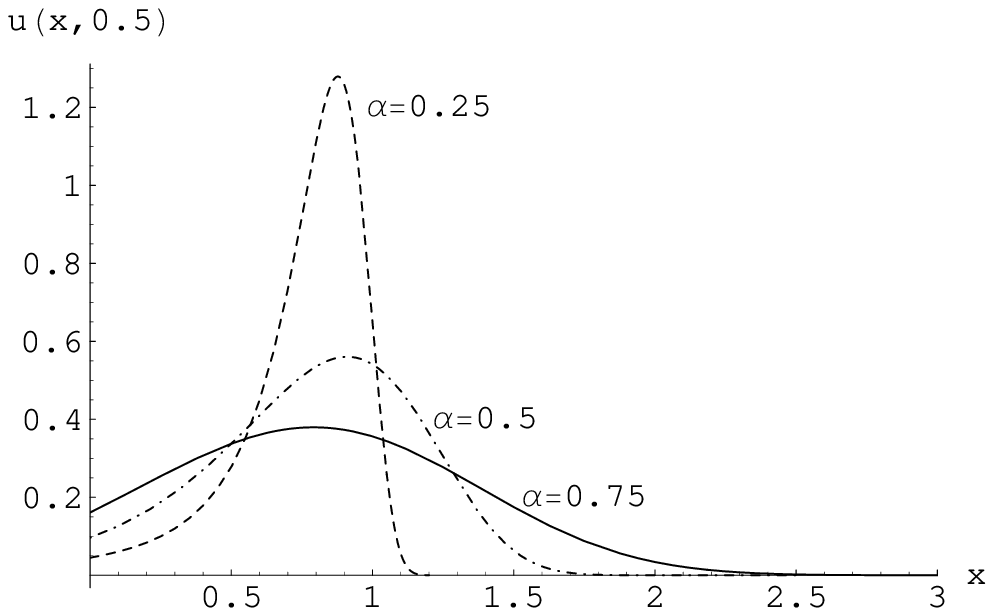}
 \includegraphics[scale=0.5]{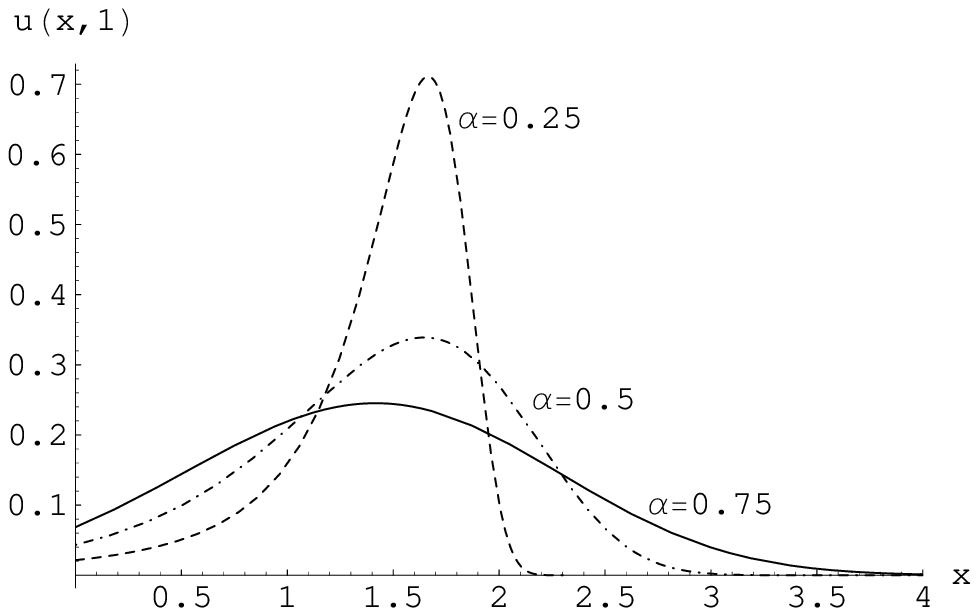}
 \includegraphics[scale=0.5]{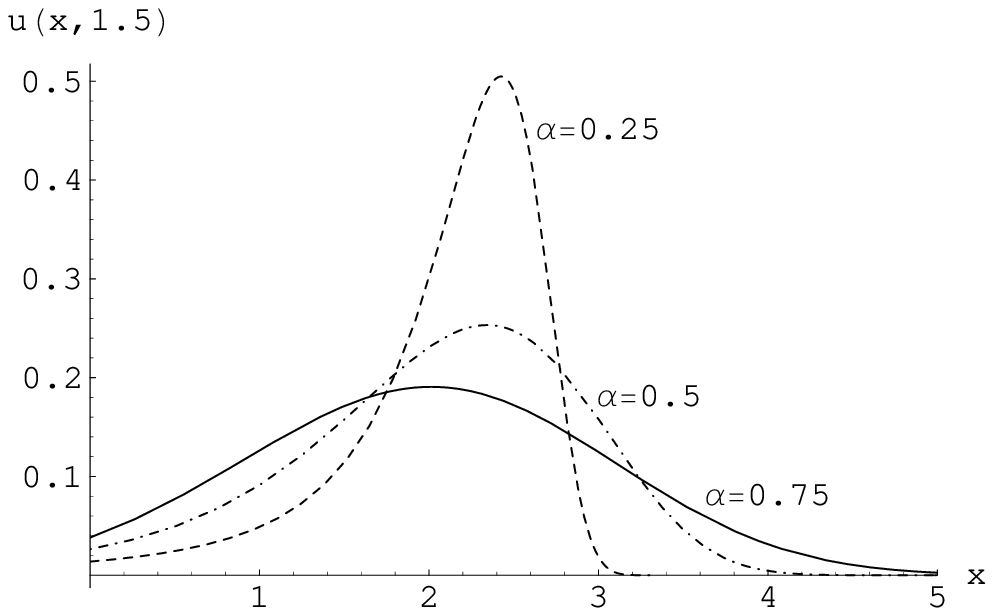}
 \end{array}$
 \end{center}
 \caption{Solution $u(x,t)$, $x\in (0,5)$, $t\in \{0.5,1,1.5\}$,
 $\alpha \in \{0.25,0.5,0.75\}$} \label{fig-4}
 \end{figure}

 Figure \ref{fig-4} presents plots of $u(x,t)$, $x\in (0,3)$, $\alpha \in \{0.25,0.5,0.75\}$,
 given by (\ref{eq:num sol}), for different time moments, i.e. for $t\in
 \{0.5,1,1.5\}$. Parameter $\tau$ and line styles are the same as in figure \ref{fig-3}.

 From figure \ref{fig-4} one can see that at the same time instance, as
 the value of $\alpha$ increases, the maximum value of $u$
 decreases, while the width of the maximum increases. This is the
 consequence of the fact that the constitutive equation of the fractional type in
 (\ref{eq:system}) describes media that tends to be elastic as $\alpha$ tends to zero. Therefore, the
 fundamental solution (\ref{eq:num sol}) should tend to the Dirac distribution as
 $\alpha \rightarrow 0$. Also, media tends to be viscoelastic as $\alpha$ in (\ref{eq:system}) tends to
 $1$. Therefore, the fundamental solution (\ref{eq:num sol}) should tend to
 the solution of wave equation for viscoelastic media as $\alpha \rightarrow 1$.

\begin{figure}[htb]
 \begin{center}
 \includegraphics[width=8cm]{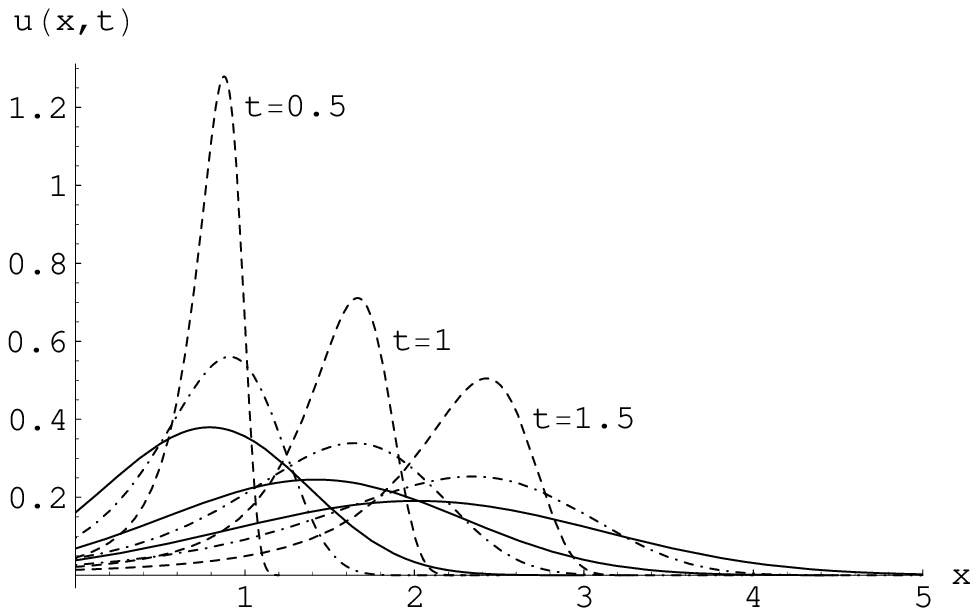}
 \end{center}
 \caption{Solution $u(x,t)$, $x\in (0,5)$, $t\in \{0.5,1,1.5\}$,
 $\alpha \in \{0.25,0.5,0.75\}$} \label{fig-5}
 \end{figure}

 Figure \ref{fig-5} presents plots of $u$ for both different time
 moments and values of $\alpha$. Parameter $\tau$ and line styles are the same as in
 figure \ref{fig-3}.

 \subsection*{Acknowledgement}

 This work is supported by Projects $144016$ and $144019$ of the
 Serbian Ministry of Science and START-project Y-237 of the
 Austrian Science Fund.


 \end{document}